\documentclass[11pt]{article}%
\usepackage{amsthm,amsmath,amsfonts}
\usepackage{amsmath}
\usepackage{amsfonts}
\usepackage{amssymb}
\usepackage{graphicx}%
\newtheorem{theorem}{Theorem}
\newtheorem{lemma}{Lemma}

\renewcommand{\d}{\mathrm{d}}

\begin{document}

\title{The limit distribution of ratios of jumps and sums of jumps of subordinators}

\author{P\'eter Kevei \\
MTA-SZTE Analysis and Stochastics Research Group \\
Bolyai Institute, Aradi v\'ertan\'uk tere 1, 6720 Szeged, Hungary \\
e-mail: \texttt{kevei@math.u-szeged.hu}
\and David M. Mason \\
Department of Applied Economics and Statistics\\University of Delaware \\
213 Townsend Hall, Newark, DE 19716, USA\\
e-mail: \texttt{davidm@udel.edu}
}

\maketitle

\begin{abstract}
Let $V_{t}$ be a driftless subordinator, and let denote $m_{t}^{(1)} \geq
m_{t}^{(2)} \geq\ldots$ its jump sequence on interval $[0,t]$. Put
$V_{t}^{(k)} = V_{t} - m_{t}^{(1)} - \ldots- m_{t}^{(k)}$ for the $k$-trimmed
subordinator. In this note we characterize under what conditions the limiting
distribution of the ratios $V_{t}^{(k)} / m_{t}^{(k+1)}$ and $m_{t}^{(k+1)} /
m_{t}^{(k)}$ exist, as $t \downarrow0$ or $t \to\infty$.
\smallskip 

\end{abstract}

\section{Introduction and results}

Let $V_{t}$, $t\geq0$, be a subordinator with L\'{e}vy measure $\Lambda$ and
drift 0. Its Laplace transform is given by
\[
\mathbf{E}\mathrm{e}^{-\lambda V_{t}}=\exp\left\{  -t\int_{0}^{\infty}\left(
1-\mathrm{e}^{-\lambda v}\right)  \Lambda(\mathrm{d}v)\right\}  ,
\]
where the L\'{e}vy measure $\Lambda$ satisfies
\begin{equation}
\int_{0}^{\infty}\min\{1,x\}\Lambda(\mathrm{d}x)<\infty.\label{min}%
\end{equation}
Put $\overline{\Lambda}(x)=\Lambda(\left(  x,\infty\right)  )$. Then
$\overline{\Lambda}(x)$ is nonincreasing and right continuous on $\left(
0,\infty\right)  $. When $t\downarrow0$ we also assume that $\overline
{\Lambda}(0+)=\infty$, which is necessary and sufficient to assure that there
is an infinite number of jumps up to time $t$, for any $t>0$.

Denote $m_{t}^{(1)}\geq m_{t}^{(2)}\geq\ldots$ the ordered jumps of $V_{s}$ up
to time $t$, and for $k\geq0$ consider the trimmed subordinator
\[
V_{t}^{(k)}=V_{t}-\sum_{j=1}^{k}m_{t}^{(j)}.
\]
We investigate the asymptotic distribution of jump sizes as $t\downarrow0$ and
$t\rightarrow\infty$. Specifically, we shall determine a necessary and
sufficient condition in terms of the L\'{e}vy measure $\Lambda$ for the
convergence in distribution of the ratios $V_{t}^{(k)}/m_{t}^{(k+1)}$ and
$m_{t}^{(k+1)}/m_{t}^{(k)}$. Observe in this notation that $V_{t}^{(0)}=V_{t}$
is the subordinator and $m_{t}^{(1)}$ is the largest jump.

An extended random variable $W$ can take the value $\infty$ with positive
probability, in which case $W$ has a defective distribution function $F$,
meaning that $F\left(  \infty\right)  <1$. We shall call an extended random
variable proper, if it is finite a.s. In this case its $F$ is a probability
distribution, i.e. $\ F\left(  \infty\right)  =1$. Here we are using the
language of the definition given on p. 127 of Feller \cite{F}.

\begin{theorem}
\label{th:jump-proc} For any choice of $k\geq0$ the ratio $V_{t}^{(k)}%
/m_{t}^{(k+1)}$ converges in distribution to an extended random variable
$W_{k}$ as $t\downarrow0$ ($t\rightarrow\infty$) if and only if one of the
following holds:

\begin{itemize}
\item[(i)] $\overline\Lambda$ is regularly varying at 0 ($\infty$) with
parameter $-\alpha$, $\alpha\in(0,1)$, in which case $W_{k}$ is a proper
random variable with Laplace transform
\begin{equation}
\label{eq:Ltrfo-limit}g_{k}(\lambda) = \frac{\mathrm{e}^{-\lambda}} {\left[  1
+ \alpha\int_{0}^{1} \left(  1 - \mathrm{e}^{-\lambda y} \right)  y^{-\alpha-
1} \d y \right]  ^{k+1}};
\end{equation}

\item[(ii)] $\overline\Lambda$ is slowly varying at 0 ($\infty$), in which
case $W_{k} = 1$ a.s.;

\item[(iii)] the condition
\begin{equation}
\frac{x\overline{\Lambda}(x)}{\int_{0}^{x}u\Lambda(\mathrm{d}u)}%
\rightarrow0\quad\text{as }x\downarrow0\ (x\rightarrow\infty
)\label{eq:infty-cond}%
\end{equation}
holds, in which case $V_{t}^{(k)}/m_{t}^{(k+1)}\overset{\mathbf{P}%
}{\longrightarrow}\infty$, that is $W_{k}=\infty$ a.s.
\end{itemize}
\end{theorem}

Note that Theorem \ref{th:jump-proc} says that the situation $0<\mathbf{P}%
\{W_{k}=\infty\}<1$ cannot happen.

The corresponding problem for nonnegative i.i.d. random variables was
investigated by Darling \cite{Darling} and Breiman \cite{Brei}, in the $k=0$
case. In this case Darling proved the sufficiency parts corresponding to (i)
and (ii) (Theorem 5.1 and Theorem 3.2 in \cite{Darling}), in particular the
limit $W_{0}$ has the same distribution as given by Darling in his Theorem
5.1, while Breiman proved the necessity parts corresponding to (i), (ii) and
(iii) (Theorem 3 (p. 357), Theorem 2 and Theorem 4 in \cite{Brei}). A special
case of Theorem 1 in Teugels \cite{Teug} gives the sufficiency analog of (i)
in the case of i.i.d. nonnegative sums for any $k\geq0$.

The necessary and sufficient condition in the cases (ii) and (iii), stated in
the more general setup of L\'{e}vy processes without a normal component, is
given by Buchmann, Fan and Maller \cite{BFM}.

\medskip Next we shall investigate the asymptotic distribution of the ratio of
two consecutive ordered jumps $m_{t}^{(k+1)}/m_{t}^{(k)}$, $k \geq 1$. We shall
obtain the analog for subordinators of a special case of a result that Bingham
and Teugels \cite{BT} established for i.i.d. nonnegative random variables.
This will follow from a general result on the asymptotic distribution of
ratios of the form defined for $k\geq1$ by%
\[
r_{k}\left(  t\right)  =\frac{\psi\left(  S_{k+1}/t\right)  }{\psi\left(
S_{k}/t\right)  }\text{, }t>0,
\]
where for each $k\geq1$, $S_{k}=\omega_{1}+\ldots+\omega_{k}$, with
$\omega_{1},\omega_{2},\ldots$ being i.i.d. mean $1$ exponential random
variables and $\psi$ is the nonincreasing and right continuous function
defined for $s>0$ by%
\begin{equation*}
\psi(s)=\sup\{y:\overline{\Pi}(y)>s\},\label{psi}%
\end{equation*}
with $\Pi$ being a positive measure on $\left(  0,\infty\right)  $ such that
$\overline{\Pi}(x)=\Pi\left(  \left(  x,\infty\right)  \right)  $
$\rightarrow0$, as $x\rightarrow\infty$. Note that we do not require $\Pi$ to
be a L\'{e}vy measure. Also whenever we consider the asymptotic distribution
of $r_{k}\left(  t\right)  $ as $t\downarrow0$ we shall assume that
$\overline{\Pi}(0+)=\infty$.

We call a function $f$ \textit{rapidly varying at 0} with index $-\infty$,
$f \in \mathrm{RV}_0(-\infty)$, if
\[
\lim_{x \downarrow 0} \frac{f(\lambda x)}{f(x)} =
\begin{cases}
0, & \text{for } \lambda > 1, \\
1, & \text{for } \lambda = 1, \\
\infty, & \text{for } \lambda < 1.
\end{cases}
\]
Correspondingly, a function $f$ is \textit{rapidly varying at $\infty$}
with index $-\infty$, $f \in \mathrm{RV}_\infty(-\infty)$, if the same holds
with $x \to \infty$.

\begin{theorem}
\label{th:jump-jump} For any choice of $k\geq1$ the ratio $r_{k}\left(
t\right)  $ converges in distribution as $t\downarrow0$ ($t\rightarrow\infty$)
to a random variable $Y_{k}$ if and only if one of the following holds:

\begin{itemize}
\item[(i)] $\overline{\Pi}$ is regularly varying at 0 ($\infty$) with
parameter $-\alpha\in\left(  -\infty,0\right)  $, in which case $Y_{k}$ has
the Beta$(k\alpha,1)$ distribution, i.e.
\begin{equation}
G_{k}(x)=\mathbf{P}\{Y_{k}\leq x\}=x^{k\alpha},\quad x\in\lbrack0,1];
\label{eq:Gk}%
\end{equation}

\item[(ii)] $\overline{\Pi}$ is slowly varying at 0 ($\infty$), in which case
$Y_{k}=0$ a.s.

\item[(iii)] $\overline{\Pi}$ is rapidly varying at 0 ($\infty$) with
index $- \infty$, in which case $Y_{k}=1$ a.s.
\end{itemize}
\end{theorem}

Theorem \ref{th:jump-jump} has some important applications to the asymptotic
distribution of the ratio of two consecutive ordered jumps $m_{t}%
^{(k+1)}/m_{t}^{(k)}$, $k\geq 1$, of a L\'{e}vy process. Let $X_{t}$, $t\geq0$,
be a L\'{e}vy processes whose L\'{e}vy measure $\Lambda$ is concentrated on
$\left(  0,\infty\right)  $. Here in addition to $\overline{\Lambda}\left(
x\right)  \rightarrow0$ as $x\rightarrow\infty$, we require that
\begin{equation}
\int_{0}^{\infty}\min\{1,x^{2}\}\Lambda(\mathrm{d}x)<\infty.\label{min2}%
\end{equation}
In this setup one has the distributional representation for $k\geq1$
\begin{equation}
\left(  m_{t}^{(k)},m_{t}^{(k+1)}\right)  \overset{\mathcal{D}}{=}\left(
\varphi(S_{k}/t),\varphi(S_{k+1}/t)\right)  ,\label{dd}%
\end{equation}
with $\varphi$ defined for $s>0$ to be
\begin{equation}
\varphi(s)=\sup\{y:\overline{\Lambda}(y)>s\}.\label{phi}%
\end{equation}
It is readily checked that $\varphi$ is nonincreasing and right continuous.
Moreover, whenever $\Lambda$ is the L\'{e}vy measure of a subordinator
$V_{t},$ condition (\ref{min}) holds, which is equivalent to
\begin{equation}
\int_{\delta}^{\infty}\varphi(s)\mathrm{d}s<\infty\text{, for any }%
\delta>0.\label{finite}%
\end{equation}
The distributional representation in (\ref{dd}) follows from Proposition 1 in
Kevei and Mason \cite{KM}. See the proof of Theorem \ref{th:jump-jump} below,
while for general spectrally positive L\'{e}vy processes it can be deduced
using the same methods that Maller and Mason \cite{MM} derived the
distributional representation for a L\'{e}vy process given in their
Proposition 5.7.

When applying Theorem \ref{th:jump-jump} to the asymptotic distribution of
consecutive ordered jumps at $0$ or $\infty$ of a L\'{e}vy processes $X_{t}$
whose L\'{e}vy measure $\Lambda$ is concentrated on $\left(  0,\infty\right)
$, we have to keep in mind that (\ref{min2}) must always hold and (\ref{min})
must be satisfied whenever $X_{t}$ is a subordinator. For instance in the case
of a subordinator $V_{t}$, whenever $m_{t}^{(k+1)}/m_{t}^{(k)}$ converges in
distribution to a random variable $Y_{k}$ as $t\downarrow0$, Theorem
\ref{th:jump-jump} says that $\overline{\Lambda}$ is regularly varying at $0$.
Further since (\ref{min}) must hold, the parameter $-\alpha$ is necessarily be
in $\left[  -1,0\right]  $, while there is no such restriction when
considering convergence in distribution as $t\rightarrow$ $\infty$. 

In the special case when $V_{t}$ is an $\alpha$-stable subordinator,
$\alpha\in(0,1)$, and $m^{(1)}>m^{(2)}>\ldots$ is its jump sequence on
$[0,1]$, then $(m^{(1)}/V_{1},m^{(2)}/V_{1},\ldots)$ has the
Poisson--Dirichlet law with parameter $(\alpha,0)$ (\textrm{PD}$(\alpha,0)$).
See Bertoin \cite{bertoin} p. 90. The ratio of the $(k+1)^{\text{th}}$ and
$k^{\text{th}}$ element of a vector, which has the PD$(\alpha,0)$ law, has the
Beta$(k\alpha,1)$ distribution (Proposition 2.6 in \cite{bertoin}). 

\section{Proofs}

In the proofs we only consider the case when $t\downarrow0$, as the
$t\rightarrow\infty$ case is nearly identical.

\subsection{Proof of Theorem \ref{th:jump-proc}}

First we calculate the Laplace exponent of the ratio using the notation
$\varphi$ defined in (\ref{phi}). We see by the nonincreasing version of the
change of variables formula stated in (4.9) Proposition of Revuz and Yor
\cite{ry}, which is given in Lemma 1 in \cite{KM},%

\begin{equation*}%
\begin{split}
\mathbf{E}\mathrm{e}^{-\lambda V_{t}}  &  =\exp\left\{  -t\int_{0}^{\infty
}\left(  1-\mathrm{e}^{-\lambda v}\right)  \Lambda(\mathrm{d}v)\right\} \\
&  =\exp\left\{  -t\int_{0}^{\infty}\left(  1-\mathrm{e}^{-\lambda\varphi
(x)}\right)  \mathrm{d}x\right\}  .
\end{split}
\label{eq:char2form}%
\end{equation*}

The key ingredient of our proofs is a distributional representation of the
subordinator $V_{t}$ given in Kevei and Mason (Proposition 1 in \cite{KM}),
which follows from a general representation by Rosi\'{n}ski \cite{Rosinski}.
It states that for $t>0$
\begin{equation}
V_{t}\overset{\mathcal{D}}{=}\sum_{i=1}^{\infty}\varphi\left(  \frac{S_{i}}%
{t}\right)  . \label{eq:repr}%
\end{equation}
From the proof of this result it is clear that $\varphi(S_{i}/t)$
corresponds to $m_{t}^{(i)}$, for $i \geq 1$. Therefore
\[
\frac{V_{t}^{(k)}}{m_{t}^{(k+1)}}\overset{\mathcal{D}}{=}\frac{\sum
_{i=k+1}^{\infty}\varphi(S_{i}/t)}{\varphi(S_{k+1}/t)}.
\]
Conditioning on $S_{k+1}=s$ and using the independence we can write
\[%
\begin{split}
\sum_{i=k+2}^{\infty}\varphi(S_{i}/t)  &  =\sum_{i=k+2}^{\infty}\varphi\left(
\frac{s}{t}+\frac{S_{i}-s}{t}\right) \\
&  \overset{\mathcal{D}}{=}\sum_{i=1}^{\infty}\varphi\left(  \frac{s}{t}%
+\frac{S_{i}}{t}\right) \\
&  =\sum_{i=1}^{\infty}\varphi_{s/t}\left(  S_{i}/t\right)  ,
\end{split}
\]
where $\varphi_{s}(x)=\varphi(s+x)$. Note that the latter sum has the same
form as in (\ref{eq:repr}), therefore it is equal in distribution to a
subordinator $V^{(s/t)}(t)$ with Laplace transform
\begin{equation}%
\begin{split}
\mathbf{E}\mathrm{e}^{-\lambda V_{t}^{(s/t)}}  &  =\exp\left\{  -t\int
_{0}^{\infty}\left(  1-\mathrm{e}^{-\lambda\varphi_{s/t}(x)}\right)
\mathrm{d}x\right\} \\
&  =\exp\left\{  -t\int_{s/t}^{\infty}(1-\mathrm{e}^{-\lambda\varphi
(x)})\mathrm{d}x\right\}  .
\end{split}
\label{eq:trunc-sub}%
\end{equation}

Now we can compute the Laplace transform of the ratio $V_{t}^{(k)}%
/m_{t}^{(k+1)}$. Since $S_{k+1}$ has Gamma$(k+1,1)$ distribution, the law of
total probability and (\ref{eq:trunc-sub}) give
\begin{equation}%
\begin{split}
\mathbf{E}\mathrm{e}^{-\lambda\frac{V_{t}^{(k)}}{m_{t}^{(k+1)}}} &
=\mathbf{E}\mathrm{e}^{-\lambda\frac{\sum_{i=k+1}^{\infty}\varphi(S_{i}%
/t)}{\varphi(S_{k+1}/t)}}\\
&  =\int_{0}^{\infty}\frac{s^{k}}{k!}\mathrm{e}^{-s}\left[  \mathrm{e}%
^{-\lambda}\mathbf{E}\mathrm{e}^{-\frac{\lambda}{\varphi(s/t)}\sum
_{i=1}^{\infty}\varphi_{s/t}(S_{i}/t)}\right]  \mathrm{d}s\\
&  =\mathrm{e}^{-\lambda}\int_{0}^{\infty}\frac{s^{k}}{k!}\mathrm{e}^{-s}%
\exp\left\{  -t\int_{s/t}^{\infty}\left[  1-\mathrm{e}^{-\frac{\lambda
}{\varphi(s/t)}\varphi(x)}\right]  \mathrm{d}x\right\}  \mathrm{d}s\\
&  =\frac{t^{k+1}}{k!}\mathrm{e}^{-\lambda}\int_{0}^{\infty}u^{k}\exp\left\{
-t\left(  u+\int_{u}^{\infty}\left[  1-\mathrm{e}^{-\lambda\frac{\varphi
(x)}{\varphi(u)}}\right]  \mathrm{d}x\right)  \right\}  \mathrm{d}u\\
&  =\frac{t^{k+1}}{k!}\mathrm{e}^{-\lambda}\int_{0}^{\infty}u^{k}%
\mathrm{e}^{-t\Psi(u,\lambda)}\mathrm{d}u,
\end{split}
\label{eq:char-funct}%
\end{equation}
where
\begin{equation}
\Psi(u,\lambda)=u+\int_{u}^{\infty}[1-\mathrm{e}^{-\lambda\frac{\varphi
(x)}{\varphi(u)}}]\mathrm{d}x.\label{eq:psi-def}%
\end{equation}
Since $\varphi$ is right continuous on $\left(  0,\infty\right)  $,
$\Psi(\cdot,\lambda)$ is also right continuous on $\left(  0,\infty\right)  $.
Further a short calculation shows that this function is strictly increasing
for any $\lambda>0$, moreover for $u_{1}>u_{2}$
\[
\Psi(u_{1},\lambda)-\Psi(u_{2},\lambda)\geq\mathrm{e}^{-\lambda}(u_{1}-u_{2}).
\]
Clearly $\Psi(0,\lambda)=0$ and $\Psi(\infty,\lambda)=\infty$. Therefore
\[
\Psi_{k}\left(  \cdot,\lambda\right)  :=\Psi \left(\left(  (k+1)\cdot\right)
^{1/(k+1)},\lambda \right)
\]
has a right continuous increasing inverse function given by
\[
Q_{\lambda}(s)=\inf\left\{  v:\Psi_{k}\left(  v,\lambda\right)  >s\right\}
\text{, for }s\geq0\text{,}%
\]
such that $Q_{\lambda}(0)=0$ and $\lim_{x\rightarrow\infty}Q_{\lambda
}(x)=\infty$. (For the right continuity part see (4.8) Lemma in Revuz and Yor
\cite{ry}.)

\medskip

\noindent\textbf{Necessity.} Assuming that $V_{t}^{(k)}/m_{t}^{(k+1)}$
converges in distribution as $t\rightarrow0$ to some extended random variable
$W_{k}$, we can apply Theorem 2a on p. 210 of Feller \cite{F} to conclude that
its Laplace transform also converges, i.e.%
\[
\int_{0}^{\infty}u^{k}\mathrm{e}^{-t\Psi(u,\lambda)}\mathrm{d}u=\int
_{0}^{\infty}\mathrm{e}^{-t\Psi_{k}\left(  v,\lambda\right)  }\mathrm{d}v
\]%
\[
=\int_{0}^{\infty}\mathrm{e}^{-ty}\mathrm{d}Q_{\lambda}\left(  y\right)
\sim\frac{\mathrm{e}^{\lambda}g_{k}(\lambda)k!}{t^{k+1}}\text{, as
}t\rightarrow0\text{,}%
\]
where $g_{k}(\lambda)=\mathbf{E}\mathrm{e}^{-\lambda W_{k}}$, and $W_{k}$ can
possibly have a defective distribution, i.e. possibly $\mathbf{P}\left\{
W_{k}=\infty\right\}  >0$. (Here we used the change of variables formula given
in (4.9) Proposition in Revuz and Yor \cite{ry}.) By Karamata's Tauberian
theorem (Theorem 1.7.1 in \cite{BGT})
\[
Q_{\lambda}(y)\sim\frac{y^{k+1}}{k+1}\mathrm{e}^{\lambda}g_{k}(\lambda
),\quad\text{as }y\rightarrow\infty,
\]
and thus by Theorem 1.5.12 in \cite{BGT}%
\[
\Psi_{k}\left(  v,\lambda\right)  \sim\left(  \frac{(k+1)v}{\mathrm{e}%
^{\lambda}g_{k}(\lambda)}\right)  ^{1/(k+1)},\text{ \ \ as }v\rightarrow
\infty,
\]
and hence
\[
\Psi(u,\lambda)\sim u\left[  \mathrm{e}^{\lambda}g_{k}(\lambda)\right]
^{-\frac{1}{k+1}},\text{ \ as }u\rightarrow\infty.
\]
Substituting back into (\ref{eq:psi-def}) we obtain for any $\lambda>0$
\begin{equation}
\lim_{u\rightarrow\infty}\frac{1}{u}\int_{u}^{\infty}\left(  1-\mathrm{e}%
^{-\lambda\frac{\varphi(x)}{\varphi(u)}}\right)  \mathrm{d}x=\left[
\mathrm{e}^{\lambda}g_{k}(\lambda)\right]  ^{-\frac{1}{k+1}}%
-1.\label{eq:asy-1}%
\end{equation}

Note that the limit $W_{k}$ is $\geq1$, with probability 1, and so
$g_{k}(\lambda)\leq\mathrm{e}^{-\lambda}$. Thus for any $\lambda$
\[
\left[  \mathrm{e}^{\lambda}g_{k}(\lambda)\right]  ^{-\frac{1}{k+1}}-1\geq0.
\]

For any $x\geq0$ we have $1-\mathrm{e}^{-x}\leq x$. Therefore by
(\ref{eq:asy-1}) we obtain for any $\lambda>0$
\begin{equation}
\liminf_{u\rightarrow\infty}\frac{1}{u\varphi(u)}\int_{u}^{\infty}%
\varphi(x)\mathrm{d}x\geq\frac{1}{\lambda}\left(  \left[  \mathrm{e}^{\lambda
}g_{k}(\lambda)\right]  ^{-\frac{1}{k+1}}-1\right)  . \label{eq:liminf}%
\end{equation}

On the other hand, by monotonicity $\varphi(x)/\varphi(u)\leq1$ for $u\leq x$.
Therefore for any $1>\varepsilon>0$ there exists a $\lambda_{\varepsilon}>0$,
such that for all $0<\lambda<\lambda_{\varepsilon}$
\[
1-\mathrm{e}^{-\lambda\frac{\varphi(x)}{\varphi(u)}}\geq(1-\varepsilon
)\frac{\lambda\varphi(x)}{\varphi(u)}\text{, for }x\geq u\text{.}%
\]
Using again (\ref{eq:asy-1}) and keeping (\ref{finite}) in mind, this implies
that for such $\lambda$
\begin{equation}
\limsup_{u\rightarrow\infty}\frac{1}{u\varphi(u)}\int_{u}^{\infty}%
\varphi(x)\mathrm{d}x\leq\frac{1}{1-\varepsilon}\frac{1}{\lambda}\left(
\left[  \mathrm{e}^{\lambda}g_{k}(\lambda)\right]  ^{-\frac{1}{k+1}}-1\right)
.\label{eq:limsup}%
\end{equation}
In particular, we obtain that, whenever $g_{k}(\lambda)\not \equiv 0$
(i.e.~$\mathbf{P}\{W_{k}<\infty\}>0$)
\[
0\leq\liminf_{u\rightarrow\infty}\frac{1}{u\varphi(u)}\int_{u}^{\infty}%
\varphi(x)\mathrm{d}x\leq\limsup_{u\rightarrow\infty}\frac{1}{u\varphi(u)}%
\int_{u}^{\infty}\varphi(x)\mathrm{d}x<\infty.
\]
Note that in (\ref{eq:liminf}) the greatest lower bound is $0$ for all
$\lambda>0$ if and only if $g_{k}(\lambda)=\mathrm{e}^{-\lambda}$, in which
case $W_{k}=1$. Then the upper bound for the limsup in (\ref{eq:limsup}) is
$0$, thus
\[
\lim_{u\rightarrow\infty}\frac{1}{u\varphi(u)}\int_{u}^{\infty}\varphi
(x)\mathrm{d}x=0,
\]
which by Proposition 2.6.10 in \cite{BGT} applied to the function 
$f(x) = x \varphi(x)$ implies that
$\varphi \in \mathrm{RV}_\infty(-\infty)$, and so, by
Theorem 2.4.7 in \cite{BGT}, $\overline{\Lambda}$ is slowly varying at $0$. We
have proved that $W_{k}=1$ if and only if $\overline{\Lambda}$ is slowly
varying at $0$.

In the following we assume that $\mathbf{P}\left\{  W_{k}>1\right\}  >0$,
therefore the liminf in (\ref{eq:liminf}) is strictly positive. Let
\begin{equation*}
a=\liminf_{\lambda\downarrow0}\frac{1}{\lambda}\left(  \left[  \mathrm{e}%
^{\lambda}g_{k}(\lambda)\right]  ^{-\frac{1}{k+1}}-1\right)  \leq
\limsup_{\lambda\downarrow0}\frac{1}{\lambda}\left(  \left[  \mathrm{e}%
^{\lambda}g_{k}(\lambda)\right]  ^{-\frac{1}{k+1}}-1\right)  =b. \label{eq:ab}%
\end{equation*}
By (\ref{eq:limsup}) and (\ref{eq:liminf}), $a>0$ and $b<\infty$. Moreover
\[
b\leq\liminf_{u\rightarrow\infty}\frac{1}{u\varphi(u)}\int_{u}^{\infty}%
\varphi(x)\mathrm{d}x\leq\limsup_{u\rightarrow\infty}\frac{1}{u\varphi(u)}%
\int_{u}^{\infty}\varphi(x)\mathrm{d}x\leq a,
\]
which forces
\[
a=b=\lim_{u\rightarrow\infty}\frac{1}{u\varphi(u)}\int_{u}^{\infty}%
\varphi(x)\mathrm{d}x=\lim_{\lambda\downarrow0}\frac{1}{\lambda}\left(
\left[  \mathrm{e}^{\lambda}g_{k}(\lambda)\right]  ^{-\frac{1}{k+1}}-1\right)
.
\]
By Karamata's theorem (Theorem 1.6.1 (ii) in \cite{BGT}) we obtain that
$\varphi$ is regularly varying at infinity with parameter $-a^{-1}%
-1=:-\alpha^{-1}$, so $\Lambda$ is regularly varying with parameter $-\alpha$
at zero with $\alpha\in(0,1)$. \smallskip

Let us consider the case when $W_{k}=\infty$ a.s., that is $V_{t}^{(k)}%
/m_{t}^{(k+1)}\overset{\mathbf{P}}{\longrightarrow}\infty$. All the previous
computations are valid, with $g_{k}(\lambda)=\mathbf{E}\mathrm{e}%
^{-\lambda\infty}\equiv0$. Thus, from (\ref{eq:liminf}) we have
\begin{equation*}
\lim_{u\rightarrow\infty}\frac{1}{u\varphi(u)}\int_{u}^{\infty}\varphi
(x)\mathrm{d}x=\infty. \label{eq:infty}%
\end{equation*}
From this, through the change of variables formula we obtain
(\ref{eq:infty-cond}).

\medskip\noindent\textbf{Sufficiency and the limit.} Consider first the
special case when $\varphi(x)=x^{-\frac{1}{\alpha}}$, $\alpha\in(0,1)$. Then a
quick calculation gives
\[
\frac{1}{u}\int_{u}^{\infty}\left(  1-\mathrm{e}^{-\lambda\frac{\varphi
(x)}{\varphi(u)}}\right)  \mathrm{d}x=\alpha\int_{0}^{1}\left(  1-\mathrm{e}%
^{-\lambda y}\right)  y^{-\alpha-1}\mathrm{d}y.
\]
By formula (\ref{eq:asy-1}) for the Laplace transform of the limit we obtain
(\ref{eq:Ltrfo-limit}).

The sufficiency can be proved by standard arguments for regularly varying
functions. Using Potter bounds (Theorem 1.5.6 in \cite{BGT}) one can show
that for $\alpha \in (0,1)$
\[
\lim_{u\rightarrow\infty}\frac{1}{u}\Psi(u,\lambda)=1+\alpha\int_{0}%
^{1}\left(  1-\mathrm{e}^{-\lambda y}\right)  y^{-\alpha-1}\mathrm{d}y,
\]
from which, through formula (\ref{eq:char-funct}), the convergence readily follows.
As already mentioned, cases (ii) and (iii) are treated in \cite{BFM}.

\subsection{Proof of Theorem \ref{th:jump-jump}}

Using that $\psi(s)\leq x$ if and only if $\overline{\Pi}(x)\leq s$, for the
distribution function of the ratio we have for $x\in(0,1)$
\begin{equation}%
\begin{split}
\mathbf{P}\left\{  r_{k}\left(  t\right)  \leq x\right\}   &  =\mathbf{P}%
\left\{  \frac{\psi(S_{k+1}/t)}{\psi(S_{k}/t)}\leq x\right\} \\
&  =\int_{0}^{\infty}\frac{s^{k-1}}{(k-1)!}\mathrm{e}^{-s}\mathbf{P}\left\{
\psi\left(  \frac{s+S_{1}}{t}\right)  \leq x\psi\left(  \frac{s}{t}\right)
\right\}  \mathrm{d}s\\
&  =\int_{0}^{\infty}\frac{s^{k-1}}{(k-1)!}\mathrm{e}^{-s}\mathrm{e}%
^{-[t\overline{\Pi}(x\psi(s/t))-s]}\mathrm{d}s\\
&  =\frac{t^{k}}{(k-1)!}\int_{0}^{\infty}u^{k-1}\mathrm{e}^{-t\overline{\Pi
}(x\psi(u))}\mathrm{d}u.
\end{split}
\label{eq:df}%
\end{equation}

\noindent\textbf{Necessity.} Assume that the limit distribution function
$G_{k}$ exists. Write
\[
\frac{t^{k}}{(k-1)!}\int_{0}^{\infty}u^{k-1}\mathrm{e}^{-t\overline{\Pi}%
(x\psi(u))}\mathrm{d}u=\frac{t^{k}}{(k-1)!}\int_{0}^{\infty}\mathrm{e}%
^{-t\Phi\left(  v,x\right)  }\mathrm{d}v,
\]
where $\Phi\left(  \cdot,x\right)  =\overline{\Pi}\left(x\psi(\left(
k\cdot\right)^{1/k})\right)$. Note that for each $x\in(0,1)$ the function
$\Phi\left(  \cdot,x\right)  $ is monotone nonincreasing and right continuous,
since $\overline{\Pi}$ and $\psi$ are both monotone nonincreasing and right
continuous. Let
\[
\mathcal{G}_{k}=\left\{  x:x\text{ is a continuity point of }G_{k}\text{ in
}(0,1)\text{ such that }G_{k}(x)>0\right\}  .
\]
First assume that $\mathbf{P}\{Y_{k}<1\}>0$. Clearly we can now proceed as in
the proof of Theorem 1 to apply Karamata's Tauberian theorem (Theorem 1.7.1 in
\cite{BGT}) to give that for any$\ x\in\mathcal{G}_{k}$,
\begin{equation}
\lim_{u\rightarrow\infty}\frac{\overline{\Pi}(x\psi(u))}{u}=[G_{k}%
(x)]^{-\frac{1}{k}}.\label{eq:asy-2}%
\end{equation}

We claim that (\ref{eq:asy-2}) implies the regular variation of $\overline
{\Pi}$. When $\overline{\Pi}$ is continuous and strictly decreasing we get by
changing variables to $\psi(u)=t$, $u=\overline{\Pi}(t)$, that we have for any
$x\in\mathcal{G}_{k}$
\[
\lim_{t\downarrow0}\frac{\overline{\Pi}(tx)}{\overline{\Pi}(t)}=[G_{k}%
(x)]^{-\frac{1}{k}},
\]
which by an easy application of Proposition 1.10.5 in \cite{BGT} implies that
$\overline{\Pi}$ is regularly varying.

Note that the jumps of $\overline{\Pi}$ correspond to constant parts of $\psi
$, and vice versa. Put $\mathcal{J}=\{z:\overline{\Pi}(z-)>\overline{\Pi
}(z)\}$ for the jump points of $\overline{\Pi}$. For $z\in\mathcal{J}$ and
$y\in\left[  \overline{\Pi}(z),\overline{\Pi}(z-)\right)  $ we have
$\psi(y)=z$. Substituting into (\ref{eq:asy-2}) we have
\begin{equation}
\lim_{z\downarrow0,z\in\mathcal{J}}\frac{\overline{\Pi}(xz)}{\overline{\Pi
}(z)}=[G_{k}(x)]^{-\frac{1}{k}},\ \text{ and }\lim_{z\downarrow0,z\in
\mathcal{J}}\frac{\overline{\Pi}(xz)}{\overline{\Pi}(z-)}=[G_{k}%
(x)]^{-\frac{1}{k}}. \label{L}%
\end{equation}
To see how the second limit holds in (\ref{L}) note that for any
$0<\varepsilon<1$ and $z\in\mathcal{J}$, we have $\psi\left(  \varepsilon
\overline{\Pi}(z)+\left(  1-\varepsilon\right)  \overline{\Pi}(z-)\right)  =z$
and thus
\[
\lim_{z\downarrow0,z\in\mathcal{J}}\frac{\overline{\Pi}(xz)}{\varepsilon
\overline{\Pi}(z)+\left(  1-\varepsilon\right)  \overline{\Pi}(z-)}%
=[G_{k}(x)]^{-\frac{1}{k}}.
\]
Since $0<\varepsilon<1$ can be chosen arbitrarily close to $0$ this implies
the validity of the second limit in (\ref{L}). Therefore by choosing any
$x\in\mathcal{G}_{k}$ we get
\begin{equation}
\lim_{z\downarrow0}\frac{\overline{\Pi}(z-)}{\overline{\Pi}(z)}=1.
\label{eq:tail-eq}%
\end{equation}

Let
\[
\mathcal{A}=\{z>0:\,\overline{\Pi}(z-\varepsilon)>\overline{\Pi}(z)\ \text{
for all }z>\varepsilon>0\}.
\]
This set contains exactly those points $z$ for which $\psi(\overline{\Pi
}(z))=z$. With this notation formula (\ref{eq:asy-2}) can be written as
\begin{equation}
\lim_{z\downarrow0,z\in\mathcal{A}}\frac{\overline{\Pi}(xz)}{\overline{\Pi
}(z)}=[G_{k}(x)]^{-\frac{1}{k}},\text{ for }x\in\mathcal{G}_{k}\text{.}
\label{G}%
\end{equation}
This together with (\ref{eq:tail-eq}) will allow us to apply Proposition
1.10.5 in \cite{BGT} to conclude that $\overline{\Pi}$ is regularly varying.
We shall need the following technical lemma.\medskip

\begin{lemma}
Whenever (\ref{eq:tail-eq}) holds, there exists a strictly decreasing sequence
$z_{n}\in\mathcal{A}$ such that $z_{n}\rightarrow0$ and
\begin{equation}
\lim_{n\rightarrow\infty}\frac{\overline{\Pi}(z_{n+1})}{\overline{\Pi}(z_{n}%
)}=1. \label{one}%
\end{equation}

\end{lemma}

\textit{Proof. }Choose $z_{1}\in\mathcal{A}$ such that $\overline{\Pi}%
(z_{1})>0$, and define for each $n\geq1$
\[
z_{n+1}=\sup\left\{  z>0:\overline{\Pi}(z)>\left(  1+\frac{1}{n}\right)
\overline{\Pi}(z_{n}-)\right\}  .
\]
Notice that the sequence $\left\{  z_{n}\right\}  $ is well-defined, since
$\overline{\Pi}(0+)=\infty$ and it is decreasing. Further we have
\[
\overline{\Pi}(z_{n+1}-)\geq\left(  1+\frac{1}{n}\right)  \overline{\Pi}%
(z_{n}-)\text{ and }\overline{\Pi}(z_{n+1})\leq\left(  1+\frac{1}{n}\right)
\overline{\Pi}(z_{n}-),
\]
where the second inequality follows by right continuity of $\overline{\Pi}$.
Also note that $z_{n+1}<z_{n}$, since otherwise if $z_{n+1}=z_{n}$, then
\[
\overline{\Pi}(z_{n+1}-)=\overline{\Pi}(z_{n}-)>\left(  1+\frac{1}{n}\right)
\overline{\Pi}(z_{n}-)\text{,}%
\]
which is impossible. Observe that each $z_{n+1}$ is in $\mathcal{A}$ since by
the definition of $z_{n+1}$ for all $0<\varepsilon<$ $z_{n+1}$%
\[
\overline{\Pi}(z_{n+1}-\varepsilon)>\left(  1+\frac{1}{n}\right)
\overline{\Pi}(z_{n}-)\geq\overline{\Pi}(z_{n+1}).
\]
Clearly since $\left\{  z_{n}\right\}  $ is a decreasing and positive
sequence, $\lim_{n\rightarrow\infty}z_{n}=z^{\ast}$ exists and is $\geq0$. By
construction
\[
\overline{\Pi}(z_{n+1}-)\geq\left(  1+\frac{1}{n}\right)  \overline{\Pi}%
(z_{n}-)\text{ }\geq\prod_{k=1}^{n}\left(  1+\frac{1}{k}\right)  \overline
{\Pi}(z_{1}-).
\]
The infinite product $\prod_{n=1}^{\infty}(1+1/n)=\infty$ forces $z^{\ast}=0$.
Also by construction we have
\[
1\leq\frac{\overline{\Pi}(z_{n+1})}{\overline{\Pi}(z_{n}-)}=\frac
{\overline{\Pi}(z_{n+1})}{\overline{\Pi}(z_{n})}\left(  \frac{\overline{\Pi
}(z_{n})}{\overline{\Pi}(z_{n}-)}\right)  \leq1+\frac{1}{n}.
\]
By (\ref{eq:tail-eq}) we have
\[
\lim_{n\rightarrow\infty}\frac{\overline{\Pi}(z_{n})}{\overline{\Pi}(z_{n}%
-)}=1.
\]
Therefore we get (\ref{one}). $\sqcup\!\!\!\!\sqcap$ \medskip

According to Proposition 1.10.5 in \cite{BGT} to establish that $\overline
{\Pi}$ is regularly varying at zero it suffices to produce $\lambda_{1}$ and
$\lambda_{2}$ in $(0,1)$ such that for $i=1,2$
\[
\frac{\overline{\Pi}(\lambda_{i}z_{n})}{\overline{\Pi}(z_{n})}\rightarrow
d_{i}\in\left(  0,\infty\right)  \text{, as }n\rightarrow\infty\text{,}%
\]
where $\left(  \log\lambda_{1}\right)  /\left(  \log\lambda_{2}\right)  $ is
finite and irrational. This can clearly be done using (\ref{G}) and
$\mathbf{P}\{Y_{k}<1\}>0$. Necessarily $\overline{\Pi}$ has index of regular
variation parameter $-\alpha\in(-\infty,0]$. For $\alpha\in(0,\infty)$ the
limiting distribution function has the form (\ref{eq:Gk}). In the case
$\alpha=0$, $\overline{\Pi}$ is slowly varying at 0 and we get that
$G_{k}(x)=1$ for $x\in(0,1)$, i.e.~$W_{k}=0$ a.s.

Now consider the case when $\mathbf{P}\{Y_{k}=1\}=1$, i.e.~$G_{k}(x)=0$ for
any $x\in(0,1)$. We once more use Theorem 1.7.1 in \cite{BGT} with $c=0$ this
time, and as an analog of (\ref{eq:asy-2}) we obtain
\[
\lim_{u\rightarrow\infty}\frac{\overline{\Lambda}(x\psi(u))}{u}=\infty.
\]
This readily implies that
\[
\lim_{z\downarrow0,z\in\mathcal{A}}\frac{\overline{\Lambda}(xz)}%
{\overline{\Lambda}(z)}=\infty,
\]
from which $\overline{\Lambda} \in \mathrm{RV}_0(-\infty)$ follows along the same
lines as before.\smallskip

\noindent\textbf{Sufficiency.}
Assume that $\overline{\Pi}$ is regularly
varying at $0$ with index $-\alpha \in (-\infty, 0)$. Then its asymptotic inverse function
$\psi$ is regularly varying at $\infty$ with index $-1/\alpha$, therefore
simply
\[
r_k(t) = \frac{\psi(S_{k+1}/t)}{\psi(S_{k}/t)} \rightarrow
\left( \frac{S_{k}}{S_{k+1}} \right)^{1/\alpha} \quad  \text{a.s., as } t \downarrow 0,
\]
which has the distribution $G_k$ in (\ref{eq:Gk}).
Assume now that $\overline{\Pi}$ is slowly varying at 0. Then $\psi \in \mathrm{RV}_\infty(-\infty)$,
therefore
\[
r_k(t) = \frac{\psi(S_{k+1}/t)}{\psi(S_{k}/t)} \rightarrow
0 \quad  \text{a.s., as } t \downarrow 0. 
\]
Finally, if $\overline{\Pi} \in \mathrm{RV}_0(-\infty)$ then $\psi$ is slowly varying
at infinity, so
\[
r_k(t) = \frac{\psi(S_{k+1}/t)}{\psi(S_{k}/t)} \rightarrow
1 \quad  \text{a.s., as } t \downarrow 0, 
\]
and the theorem is completely proved.

\medskip

\noindent \textbf{Acknowledgement.}
PK was supported by the Hungarian
Scientific Research Fund OTKA PD106181 and by the European Union and
co-funded by the European Social Fund under the project
`Telemedicine-focused research activities on the field of Mathematics,
Informatics and Medical sciences' of project number T\'{A}MOP-4.2.2.A-11/1/KONV-2012-0073.
DM thanks the Bolyai Institute for their hospitality while this paper was being written.

\end{document}